\documentclass{kapproc}
\usepackage[T2A]{fontenc}
\usepackage[cp866]{inputenc}
\usepackage[russian]{babel}
\usepackage{amsmath}
\usepackage{amssymb}

\newtheorem{theo}{Theorem}
\newtheorem{prop}{Proposition}
\newtheorem{lemma}{Lemma}
\newtheorem{corol}{Corollary}[theo]
\def\bydef{\stackrel{\mathrm{def}}{=}}
\def\N{\mathbb{N}}
\def\R{\mathbf{R}}
\def\C{\mathbb{C}}

\def\PN{\mathbb{P}}
\def\CN{\mathbb{C}}
\def\LS{\mathrm{L}}
\def\LN{\mathbb{L}}
\def\HS{\mathrm{H}}
\def\HN{\mathbb{H}}
\def\TS{\mathrm{T}}
\def\TN{\mathbb{T}}
\def\F{{\cal F}}
\def\AH{{\cal H}_{\theta}}
\def\AL{{\cal L}_{\theta}^{\nu}}
\def\AJ{{\cal J}_{\theta}^{\alpha,\beta}}
\def\AC{{\cal G}_{\theta}^{\lambda}}
\def\HK{H\!K}
\def\LK{L\!K}
\def\JK{J\!K}
\def\GK{G\!K}
\def\bydef{\stackrel{\mathrm{def}}{=}}

\normallatexbib

\begin{document}
\articletitle{Square summability with geometric weight for classical
orthogonal expansions}

\author{Dmitrii Karp}
\affil{Institute of Applied Mathematics\\7 Radio Street, Vladivostok,
690041, Russia} \email{dmkrp@yandex.ru}

\begin{abstract}
Let $f_k$ be the $k$-th Fourier coefficient of a function $f$ in
terms of the orthonormal Hermite, Laguerre or Jacobi polynomials. We
give necessary and sufficient conditions on $f$ for the inequality
$\sum_{k}|f_k|^2\theta^k<\infty$ to hold with $\theta>1$. As a
by-product new orthogonality relations for the Hermite and Laguerre
polynomials are found. The basic machinery for the proofs is provided
by the theory of reproducing kernel Hilbert spaces.
\end{abstract}



\paragraph{1. Introduction}
The goal of this paper is to find necessary and sufficient conditions
to be imposed on a function $f$ for its Fourier coefficients in terms
of classical orthogonal polynomials to satisfy the inequality
\begin{equation}\label{eq:sum}
\sum_{k=0}^{\infty}|f_k|^2\theta^k<\infty
\end{equation}
with $\theta>1$.  So, we have three problems corresponding to the
following three definitions of $f_k$:
\begin{equation}\label{eq:HermiteCoef}
f_k=\int\limits_{-\infty}^{\infty}\!\!f(x)\HN_k(x)e^{-x^2}dx,
\end{equation}
\begin{equation}\label{eq:LaguerreCoef}
f_k=\int\limits_{0}^{\infty}\!\!f(x)\LN^{\!\nu}_k(x)x^{\nu}e^{-x}dx,
\end{equation}
and
\begin{equation}\label{eq:JacobiCoef}
f_k=\int\limits_{-1}^{1}\!\!f(x)\PN^{\alpha,\beta}_{k}(x)(1-x)^{\alpha}(1+x)^{\beta}dx.
\end{equation}
Here $\HN_k$, $\LN^{\!\nu}_{k}$ and $\PN^{\alpha,\beta}_{k}$ is the
$k$-th orthonormal polynomial of Hermite, Laguerre and Jacobi,
respectively \cite{Szego91}.  For the sake of convenience we use
orthonormal instead of standardly normalized versions of the
classical polynomials. In each case $f$ is defined on the interval of
orthogonality of the corresponding system of polynomials.  We will
use $\varphi_k$ as a generic notation for either of the three types
of polynomials.

Classes of functions with rapidly decreasing Fourier coefficients in
classical orthogonal polynomials have been extensively studied. We
only mention a few contributions without any attempt to make a
survey.  The series of papers \cite{Hille39}-\cite{Hille80} by
E.~Hille is devoted to the Hermite expansions.  Among other things
Hille studied expansions with $f_k$ vanishing as quick as
$\exp(-\tau(2k+1)^{1/2})$. The functions possessing such coefficients
are holomorphic in the strip $|\Im{z}|<\tau$. Hille provided an exact
description of the linear vector space with compact convergence
topology formed by these functions where the set
$\{\HN_k:k\in\N\cup{0}\}$ is a basis. Its members are characterized
by a suitable growth condition. In addition, he studied the
convergence on and analytic continuation through the boundary of the
strip.

The convergence domain of the Laguerre series with
$f_k\sim\exp(-\tau{k}^{1/2})$ is the interior of the parabola
$\Re(-z)^{1/2}=\tau/2$.  Rusev in~\cite{Rusev84} described the linear
vector space with compact convergence topology formed by  functions
holomorphic in $\{z:\Re(-z)^{1/2}<\tau/2\}$ where the set
$\{\LN_k^{\nu}:k\in\N\cup{0}\}$ is a basis, and Boyadjiev
in~\cite{Boyad92} studied the behavior of the Laguerre series on the
boundary of the convergence domain.

If coefficients (\ref{eq:HermiteCoef}) or (\ref{eq:LaguerreCoef})
decrease as fast as $\exp(-\tau{k}^{\eta})$ with $\eta>1/2$ the
function $f$ is entire.  The spaces comprising such functions for
Hermite expansions have been characterized by Janssen and van
Eijndhoven in~\cite{JanEijnd90}.  If the Fourier-Hermite coefficients
$f_k$ decline quicker than any geometric progression, the suitable
characterization is provided by Berezanskij and Kondratiev
in~\cite{BerKondr95} (see Corollary~\ref{corol:Herm} below).  Byun
was the first to study the Hermite expansions with condition
(\ref{eq:sum}) - see remark after Theorem~\ref{theo:Herm}. For the
Laguerre expansions with $\limsup_{k}|f_k|^{1/k}<1$ Zayed related the
singularities of the Borel transform of $f$ with those of
$F(z)=\sum_kf_kz^k$ in~\cite{Zayed81}.

For the Legendre expansions with $\limsup_{k}|f_k|^{1/k}=\zeta<1$
Nehari relates the singularities of $f$ on the boundary of
convergence domain $\{z:|z+1|+|z-1|<\zeta+\zeta^{-1}\}$ to those of
$F(z)=\sum_kf_kz^k$ in~\cite{Nehari56}.  Gilbert in~\cite{Gilbert64}
and Gilbert and Howard in~\cite{GilbHoward66} generalized the results
of Nehari to the Gegenbauer expansions and to expansions in
eigenfunctions of a Sturm-Liouville operator.

Functions satisfying (\ref{eq:sum}) apparently form a proper subclass
of functions with
$\limsup_{k}|f_k|^{\frac{1}{k}}\leq\theta^{-\frac{1}{2}}$.  On the
other hand the condition (\ref{eq:sum}) itself cannot be expressed in
terms of asymptotics of $f_k$.  Consequently, our criteria for the
validity of (\ref{eq:sum}) are of different character from those
contained in the above references.  Although they also describe the
growth of $f$ for the Hermite and Laguerre expansions and its
boundary behaviour for the Jacobi expansions, our growth conditions
are given in terms of existence of certain weighted area integrals of
$f$ and cannot be expressed by an estimate of the modulus, while the
restriction on the boundary behaviour is given on the whole boundary
and not in terms of analysis of individual singularities.

\paragraph{2. Preliminaries} Throughout the paper the following standard notation will be
used: $\N$, $\R$, $\R^{+}$ and $\C$ will denote positive integers,
real numbers, positive real numbers and the finite complex plane,
respectively.  Since coefficients
(\ref{eq:HermiteCoef})-(\ref{eq:JacobiCoef}) do not change if we
modify $f$ on a set of Lebesgue measure zero, all our statements
about the properties of $f$ should be understood to hold almost
everywhere.  If, for instance, we say that $f$ is the restriction of
a holomorphic function to (some part of) the real axis, it means that
$f$ is allowed to differ from such restriction on a set of zero
measure.

All proofs in the paper hinge on the theory of reproducing kernel
Hilbert spaces (RKHS), so for convenience we briefly outline the
basic facts of the theory we will make use of.

For a Hilbert space $H$ comprising complex-valued functions on a set
$E$, the reproducing kernel $K(p,q): E\times{E}\to\C$ is a function
that belongs to $H$ as a function of $p$ for every fixed $q\in{E}$
and possesses the reproducing property
\[
f(q)=(f,K(\cdot,q))_{H}
\]
for every $f\in{H}$ and for any $q\in{E}$.  If a Hilbert space admits
the reproducing kernel then this kernel is unique and positive
definite on $E\times{E}$:
\begin{equation}\label{eq:KPositivity}
\sum\limits_{i,j=1}^{n}K(p_i,p_j)c_i\overline{c_j}\geq 0
\end{equation}
for an arbitrary finite complex sequence $\{c_i\}$ and any points
$p_i\in{E}$. The theorem of Moore and Aronsjain \cite{Aron50} states
that the converse is also true: every positive definite kernel $K$ on
$E\times{E}$ uniquely determines a Hilbert space $H$ admitting $K$ as
its reproducing kernel.  This fact justifies the notation $H_K$ for
the Hilbert space $H$ induced by the kernel $K$.  The following
propositions can be found in \cite{Aron50,Saitoh97}.
\begin{prop}\label{prop:Ks}
If $H_K$ is a Hilbert space of functions $E\to\C$ and $s$ is an
arbitrary non-vanishing function on $E$, then
\begin{equation}\label{eq:Ks}
K_s(p,q)=s(p)\overline{s(q)}K(p,q)
\end{equation}
is the reproducing kernel of the Hilbert space $H_{K_s}$ comprising
all functions on $E$ expressible in the form $f_s(p)=s(p)f(p)$ with
$f{\in}H_K$ and equipped with inner product
\begin{equation}\label{eq:InnerKs}
(f_s,g_s)_{H_{K_s}}=\left(\frac{f_s}{s},\frac{g_s}{s}\right)_{H_K}.
\end{equation}
\end{prop}

\begin{prop}\label{prop:Restriction}
Let $E_1\subset E$ and $K_1$ be the restriction of a positive
definite kernel $K$ to $E_1\times{E_1}$.  Then the RKHS $H_{K_1}$
comprises all restrictions to $E_1$ of functions from $H_K$ and has
the norm given by
\begin{equation}\label{eq:RestrNorm}
\|f_1\|_{H_{K_1}}=\min\left\{\|f\|_{H_K}; f\vert_{E_1}=f_1,~f\in{H_K}\right\}.
\end{equation}
\end{prop}
\begin{prop}\label{prop:KBilin}
If RKHS $H_K$ is separable and $\{\psi_k:k\in\N\}$ is a complete
orthonormal system in $H_K$, then its reproducing kernel is expressed
by
\begin{equation}\label{eq:KBilin}
K(p,q)=\sum\limits_{k}\psi_k(p)\overline{\psi_k(q)},
\end{equation}
where the series (\ref{eq:KBilin}) converges absolutely for all
$p,q\in E$ and uniformly on every subset of $E$, where $K(q,q)$ is
bounded.
\end{prop}

The relation $K_1\ll K_2$ will mean that $K_2-K_1$ is positive
definite. This relation introduces partial ordering into the set of
positive definite kernels on $E\times{E}$.  Inclusion $H_{K_1}\subset
H_{K_2}$ and equality $H_{K_1}=H_{K_2}$ will be understood in the
set-theoretic sense, which implies, however, that the same relations
hold in the topological sense as stated in the following two
propositions.
\begin{prop}\label{prop:K1subK2}
Inclusion $H_{K_1}\subset H_{K_2}$ takes place iff $K_1\ll MK_2$ for
a constant $M>0$. In this case $M^{1/2}\|f\|_{1}\geq \|f\|_{2}$ for
all $f\in H_{K_1}$. {\em(}Here $\|f\|_{1}$ and $\|f\|_{2}$ are the
norms in  $H_{K_1}$ and $H_{K_2}$, respectively{\em)}.
\end{prop}
\begin{prop}\label{prop:K1equalsK2}
Equality $H_{K_1}=H_{K_2}$ takes place iff $mK_2\ll K_1\ll MK_2$ for
some positive constants  $m$, $M$. In this case
$m^{1/2}\|f\|_{1}\leq\|f\|_2\leq M^{1/2}\|f\|_{1}$.
\end{prop}

The kernels satisfying Proposition~\ref{prop:K1equalsK2} are said to
be equivalent which is denoted by $K_1\approx{K_2}$.

It is shown in \cite{Aron50} that the RKHS $H_K$ induced by the
kernel $K(p,q)=K_1(p,q)K_2(p,q)$ consists of all restrictions to the
diagonal of $E'=E\times E$ (i.e. the set of points of the form
$\{p,p\}$) of the elements of the tensor product $H'=H_{K_1}\otimes
H_{K_2}$. The space $H_K$ is characterized by
\begin{prop}\label{prop:TensorDiagonal}
Let $K(p,q)=K_1(p,q)K_2(p,q)$ and let $\{\psi_k:k\in\N\}$ be a
complete orthonormal set in $H_{K_2}$. Then the RKHS $H_K$ comprises
the functions of the form
\begin{equation}\label{eq:DiagonalElement}
f(p)=\sum\limits_{k=1}^{\infty}f^1_k(p)\psi_k(p),~~~f^1_k\in
H_{K_1},~~\sum\limits_{k=1}^{\infty}\|f^1_k\|_{1}^2<\infty.
\end{equation}
The norm in $H_K$ is given by
\[
\|f\|_{H_{K}}^2=\min\left\{\sum\limits_{k=1}^{\infty}\|f^1_k\|_{1}^2\right\},
\]
where the minimum is taken over all representations of $f$ in the
form {\em (\ref{eq:DiagonalElement})} and is attained on one such
representation.
\end{prop}

\paragraph{3. Results for the Hermite and Laguerre expansions}
Functions satisfying (\ref{eq:sum}) form a Hilbert space with inner
product
\begin{equation}\label{eq:InnerProduct}
(f,g)=\sum\limits_{k=0}^{\infty}f_k\overline{g_k}\theta^k.
\end{equation}
This space will be denoted by $\AH$  for the Hermite expansions and
by $\AL$ for the Laguerre expansions.  The sets
$\{\theta^{-k/2}\HN_k\}_{k\in\N\cup{0}}$ and
$\{\theta^{-k/2}\LN_k^{\!\nu}\}_{k\in\N\cup{0}}$ constitute
orthonormal bases of the spaces  $\AH$ and $\AL$, respectively. For
each space we can form the reproducing kernel according to
(\ref{eq:KBilin}):
\begin{equation}\label{eq:KBilinpolynom}
K(z,\overline{u})=\sum\limits_{k}\varphi_k(z)\overline{\varphi_k(u)}\theta^{-k}.
\end{equation}
The explicit formulae for these kernels are known to be
\cite{Bateman53}:
\begin{equation}\label{eq:BilinHermite}
\HK_{\theta}(z,\overline{u})\bydef\sum\limits_{k=0}^{\infty}\HN_{k}(z)\overline{\HN_{k}(u)}\theta^{-k}=
\frac{\theta}{\sqrt{\pi(\theta^2-1)}}\exp{\left(\frac{2z\overline{u}\theta-z^2-\overline{u}^2}{\theta^2-1}\right)}
\end{equation}
(Mehler's formula) and
\begin{equation}\label{eq:BilinLaguerre}
\LK_{\theta}^{\nu}(z,\overline{u})\bydef
\sum\limits_{k=0}^{\infty}\LN^{\nu}_{k}(z)\overline{\LN^{\nu}_{k}(u)}\theta^{-k}
=\frac{\theta^{\nu/2+1}}{\theta-1}\exp{\left(-\frac{z+\overline{u}}{\theta-1}\right)}(z\overline{u})^{-\frac{\nu}{2}}
I_\nu\left(2\frac{\sqrt{\theta{z\overline{u}}}}{\theta-1}\right)
\end{equation}
(Hardy-Hille's formula). Here $I_{\nu}$ is the modified Bessel
function. The kernels (\ref{eq:BilinHermite}) and
(\ref{eq:BilinLaguerre}) are entire functions of both $z$ and
$\overline{u}$. The space $\AH$ comprises functions on $\R$, while
the space $\AL$ comprises functions on $\R^+$.  Hence, we consider
the restrictions of the kernels (\ref{eq:BilinHermite}) and
(\ref{eq:BilinLaguerre}) to $\R$ and $\R^+$, respectively. Applying
Proposition~\ref{prop:Restriction} with $E=\C$ and $E'=\R$ or
$E'=\R^+$ we conclude that the spaces $\AH$ and $\AL$ are formed by
all restrictions to $\R$ and $\R^+$, respectively, of entire
functions from the spaces generated by the kernels
(\ref{eq:BilinHermite}) and (\ref{eq:BilinLaguerre}). We can drop
minimum in (\ref{eq:RestrNorm}) due to uniqueness of analytic
extension and, consequently, the norm induced by inner product
(\ref{eq:InnerProduct}) equals the norm in $H_{\HK_{\theta}}$ or in
$H_{\LK_{\theta}^{\nu}}$.

Next, we observe that both kernels (\ref{eq:BilinHermite}) and
(\ref{eq:BilinLaguerre}) are of the form
$s(z)\overline{s(u)}K(z\overline{u})$ with non-vanishing functions
$s_H(z)=e^{-z^2/(\theta^2-1)}$ and $s_L(z)=e^{-z/(\theta-1)}$. Hence
we are in the position to apply Proposition~\ref{prop:Ks}. In
compliance with (\ref{eq:InnerKs}) the norms in $H_{\HK_{\theta}}$
and $H_{\LK_{\theta}^{\nu}}$ are known once we have found the norms
in the spaces induced by the kernels
\begin{equation}\label{eq:K1Hermite}
\widetilde{\HK}_{\theta}(z\overline{u})=
\frac{\theta}{\sqrt{\pi(\theta^2-1)}}\exp{\left(\frac{2z\overline{u}\theta}{\theta^2-1}\right)}
\end{equation}
and
\begin{equation}\label{eq:K1Laguerre}
\widetilde{\LK}^{\nu}_{\theta}(z\overline{u})=\frac{\theta^{\nu/2+1}}{\theta-1}(z\overline{u})^{-\frac{\nu}{2}}
I_\nu\left(2\frac{\sqrt{\theta{z\overline{u}}}}{\theta-1}\right).
\end{equation}
Both of them depend on the product $z\overline{u}$ and are rotation
invariant thus.  Rotation invariance of the kernel implies radial
symmetry of the measure with respect to which the integral
representing the norm is taken.  For the kernel (\ref{eq:K1Hermite})
the measure and the space are well-known. It is the Fischer-Fock (or
the Bargmann-Fock) space $\F_\theta$ of entire functions with finite
norms
\begin{equation}\label{eq:Fock}
\|f\|_{\F_\theta}^2=
\frac{2}{\sqrt{\pi(\theta^2-1)}}\int\limits_{\C}\!\!|f(z)|^2\exp\left(-\frac{2\theta|z|^2}{\theta^2-1}\right)d\sigma,
\end{equation}
where the integration is with respect to Lebesgue's area measure. By
Proposition~\ref{prop:Ks} the final result for the Hermite expansions
now becomes straightforward.
\begin{theo}\label{theo:Herm}
Inequality (\ref{eq:sum}) with $f_k$ defined by
(\ref{eq:HermiteCoef}) holds true for all restrictions to $\R$ of the
entire functions with
\begin{equation}\label{eq:HermIntNorm}
\int\limits_{\C}|f(z)|^2
\exp\left[-2\left(\frac{(\Re{z})^2}{{\theta}+1}+\frac{(\Im{z})^2}{{\theta}-1}\right)\right]d\sigma<\infty
\end{equation}
and only for them.
\end{theo}
For inner product (\ref{eq:InnerProduct}) this leads to the
expression
\begin{equation}\label{eq:HermIntInner}
(f,g)_{\AH}=(f,g)_{H_{\HK_\theta}}=\frac{2}{\sqrt{\pi(\theta^2-1)}}\int\limits_{\C}f(z)\overline{g(z)}
\exp\left[-2\left(\frac{(\Re{z})^2}{{\theta}+1}+\frac{(\Im{z})^2}{{\theta}-1}\right)\right]d\sigma.
\end{equation}
Since polynomials $\HN_k$ are orthogonal with respect to this inner
product, we obtain the following orthogonality relation for the
standardly normalized  Hermite polynomials $\HS_k$:
\begin{equation}\label{eq:HermOrtho}
\frac{2}{\sqrt{(\theta^2-1)}}\int\limits_{\C}\HS_k(z)\overline{\HS_m(z)}
\exp\left[-\frac{2(\Re{z})^2}{{\theta}+1}-\frac{2(\Im{z})^2}{{\theta}-1}\right]d\sigma=
\delta_{k,m}\pi(2\theta)^kk!.
\end{equation}
These results for the Hermite expansions have been essentially proved
by Du-Wong Byun in \cite{Byun93}, although the emphasis in his work
is different and it seems that the orthogonality relation
(\ref{eq:HermOrtho}) has not been noticed.

The following corollary is an immediate consequence of
(\ref{eq:HermIntNorm}).
\begin{corol}\label{corol:Herm}
Inequality (\ref{eq:sum}) with $f_k$ defined by
(\ref{eq:HermiteCoef}) holds true for \textbf{all} $\theta>1$ iff $f$
is the restriction to $\R$ of an entire function $F$  satisfying
$|F(z)|\leq Ce^{\varepsilon |z|^2}$ for all $\varepsilon>0$ and a
constant $C=C(\varepsilon)$ independent of $z$.
\end{corol}
A direct proof is given in \cite{BerKondr95}.

For the Laguerre expansions the situation is a bit more complicated.
The kernel (\ref{eq:K1Laguerre}) is a particular case of a much more
general hypergeometric kernel.  The spaces generated by
hypergeometric kernels are studied in depth in \cite{Karp03}.  For
the kernel (\ref{eq:K1Laguerre}) we get the space of entire functions
with finite norms
\begin{equation}\label{eq:LagK1Norm}
\|f\|^2_{\widetilde{\LK}_{\theta}^{\nu}}=\frac{2\theta^{-\nu/2}}{\pi(\theta-1)}\int\limits_{\C}\!\!|f(z)|^2
|z|^{\nu}K_{\nu}\left(\frac{2\sqrt{\theta}|z|}{\theta-1}\right)d\sigma,
\end{equation}
where $K_\nu$ is the modified Bessel function of the second kind (or
the McDonald function).  Application of Proposition~\ref{prop:Ks}
brings us to our final result for the Laguerre expansions.
\begin{theo}\label{theo:Lag}
Inequality (\ref{eq:sum}) with $f_k$ defined by
(\ref{eq:LaguerreCoef}) holds true for all restrictions to $\R^+$ of
the entire functions with
\begin{equation}\label{eq:LagIntNorm}
\int\limits_{\C}\!\!|f(z)|^2
\exp\!\left(\frac{2\Re{z}}{\theta-1}\right)|z|^{\nu}
K_{\nu}\left(\frac{2\sqrt{\theta}|z|}{\theta-1}\right)d\sigma<\infty
\end{equation}
and only for them.
\end{theo}
The orthogonality relation for the standardly normalized Laguerre
polynomials $\LS^{\!\nu}_k$ that follows from this result is given by
\begin{equation}\label{eq:LagOrtho}
\frac{2\theta^{-\nu/2}}{\pi(\theta-1)}\int\limits_{\C}\!\!\LS^{\!\nu}_k(z)\overline{\LS^{\!\nu}_m(z)}
\exp\!\left(\frac{2\Re{z}}{\theta-1}\right)|z|^{\nu}
K_{\nu}\left(\frac{2\sqrt{\theta}|z|}{\theta-1}\right)d\sigma=\delta_{k,m}\frac{\Gamma(k+\nu+1)\theta^k}{k!}.
\end{equation}
\begin{corol}\label{corol:Lag}
Inequality (\ref{eq:sum}) with $f_k$ defined by
(\ref{eq:LaguerreCoef}) holds true for \textbf{all} $\theta>1$ iff
$f$ is the restriction to $\R^+$ of an entire function $F$ satisfying
$|F(z)|\leq Ce^{\varepsilon |z|}$ for all $\varepsilon>0$ and a
constant $C=C(\varepsilon)$ independent of $z$.
\end{corol}
This corollary can be easily derived from (\ref{eq:LagIntNorm})
with the help of asymptotic relation \cite{Bateman53}
\[
K_{\nu}(x)\simeq\sqrt{\frac{\pi}{2x}}e^{-x},~~x\to\infty.
\]

\paragraph{4. Results for the Jacobi expansions}
The space of complex-valued functions on $(-1,1)$ whose
Fourier-Jacobi coefficients (\ref{eq:JacobiCoef}) satisfy
(\ref{eq:sum}) will be denoted by $\AJ$.  Pursuing the same line
of argument as in the previous section, we form the reproducing
kernel of this space found by Baily's formula \cite{Bateman53}
\[
\JK^{\theta}_{\alpha,\beta}(z,\overline{u})\bydef
\sum\limits_{k=0}^{\infty}\PN^{\alpha,\beta}_k(z)\overline{\PN^{\alpha,\beta}_k(u)}\frac{1}{\theta^k}=
\frac{\theta^{\alpha+\beta+1}(\theta-1)}{\tau(\alpha,\beta)
(\theta+1)^{\alpha+\beta+2}}\times
\]
\begin{equation}\label{eq:BilinJacobi}
\times
F_4\left(\frac{\alpha}{2}+\frac{\beta}{2}+1,\frac{\alpha}{2}+\frac{\beta}{2}
+\frac{3}{2};\alpha+1,\beta+1;\frac{\theta(1-z)(1-\overline{u})}{(\theta+1)^2},
\frac{\theta(1+z)(1+\overline{u})}{(\theta+1)^2}\right),
\end{equation}
where
\begin{equation}\label{eq:AppelF4}
F_4(a,b;c,c';t,s)=\sum\limits_{m,n=0}^{\infty}\frac{(a)_{n+m}(b)_{n+m}}
{(c)_n(c')_m m! n!}t^m s^n
\end{equation}
is Appel's hypergeometric function and $\tau(\alpha,\beta)=
2^{\alpha+\beta+1}\Gamma(\alpha+1)\Gamma(\beta+1)/\Gamma(\alpha+\beta+2)$.
Define the ellipse $E_\theta$ by
\begin{equation}\label{eq:Etheta}
E_\theta=\{z:|z-1|+|z+1|<\theta^{1/2}+\theta^{-1/2}\}.
\end{equation}
Our first observation here is that both the series on the right hand
side and on the left hand side of (\ref{eq:BilinJacobi}) converge
absolutely and uniformly on compact subsets of $E_\theta\times
E_\theta$ (see \cite{Karp00}). The implication of the uniform
convergence is holomorphy of the kernel (\ref{eq:BilinJacobi}) in
$E_\theta\times{E_\theta}$ with respect to both variables.
Application of Proposition~\ref{prop:Restriction} with $E=E_\theta$,
$E'=(-1,1)$ leads to the assertion that the space $\AJ$ is formed by
 restrictions to the interval $(-1,1)$ of functions holomorphic in
$E_\theta$ thereby the norm in $\AJ$ equals the norm in
$H_{\!\JK^{\theta}_{\alpha,\beta}}$ due to uniqueness of analytic
continuation.

Our main result for the Jacobi expansions will be derived from its
particular case $\alpha=\beta=\lambda-1/2$ thereby the orthonormal
Jacobi polynomials reduce to the orthonormal Gegenbauer polynomials
$\CN^{\lambda}_k$. The reproducing kernel (\ref{eq:BilinJacobi})
reduces in this case to
\[
\GK^{\theta}_{\lambda}(z,\overline{u})\bydef\sum\limits_{k=0}^{\infty}\CN^{\lambda}_k(z)\overline{\CN^{\lambda}_k(u)}\frac{1}{\theta^{k}}=
\]
\begin{equation}\label{eq:BilinGigit}
=\frac{\theta^{2\lambda}(\theta^2-1)}{\tau(\lambda)(\theta^2-2{\theta}z\overline{u}+1)^{\lambda+1}}
{_{2}F_{1}}\left(\frac{\lambda+1}{2},\frac{\lambda+2}{2};\lambda+\frac{1}{2};\frac{4\theta^2(1-z^2)(1-\overline{u}^2)}
{(\theta^2-2{\theta}z\overline{u}+1)^2}\right),
\end{equation}
where ${_{2}F_{1}}$ is the Gauss hypergeometric function and
$\tau(\lambda)= \sqrt{\pi}\Gamma(\lambda+1/2)/\Gamma(\lambda+1)$. The
series on both sides of (\ref{eq:BilinGigit}) again converge
absolutely and uniformly on compact subsets of $E_\theta\times
E_\theta$. Formula  (\ref{eq:BilinGigit}) can be obtained from
(\ref{eq:BilinJacobi}) by using a reduction formula for $F_4$ and
applying a quadratic transformation to the resulting hypergeometric
function. Details are in \cite{Karp00}.

Let $\partial{E_{\theta}}$ denote the boundary of the ellipse
$E_\theta$. We introduce the weighted Szeg\"{o} space
$AL_2(\partial{E_{\theta}}; \rho)$ with continuous positive weight
$\rho(z)$ defined on $\partial{E_{\theta}}$ as the set of functions
holomorphic in $E_{\theta}$, possessing non-tangential boundary
values almost everywhere on $\partial{E_{\theta}}$ and having finite
norms
\begin{equation}\label{SzegoE}
\|f\|_{AL_2(\partial{E_{\theta}};\rho)}^2=\int\limits_{\partial{E_{\theta}}}\!\!|f(z)|^2\rho(z)|dz|<\infty.
\end{equation}
We will write $AL_2(\partial{E_{\theta}})$ for
$AL_2(\partial{E_{\theta}}; 1)$.

For $\lambda=0$, the orthonormal Gegenbauer polynomials reduce to the
orthonormal Chebyshev polynomials of the first kind $\TN_k$:
\begin{equation}\label{eq:T0defined}
\begin{array}{cc}
\CN^{0}_k(z)=\TN_k(z)=(2/\pi)^{\frac{1}{2}}\TS_k(z)=(2/\pi)^{\frac{1}{2}}\cos(k\arccos{z}),
& k\in\N,
\\[0.07in] \CN^{0}_0(z)=\TN_0(z)=(1/\pi)^{\frac{1}{2}}\TS_0(z)=(1/\pi)^{\frac{1}{2}}, &
\end{array}
\end{equation}
where $\TS_k$ is the $k$-th Chebyshev polynomial of the first kind in
standard normalization.
\begin{lemma}\label{lemma:T-ortho}
Polynomials $(\theta^k+\theta^{-k})^{-1/2}\TN_k/2$ form orthonormal
basis of the space
$AL_2(\partial{E_{\theta}};|z^2-1|^{-\frac{1}{2}})$.
\end{lemma}
\textbf{Proof.} The boundary $\partial{E_{\theta}}$ is an analytic
arc which implies the completeness of the set of all polynomials in
$AL_2(\partial{E_{\theta}};|z^2-1|^{-\frac{1}{2}})$ (see for instance
\cite{Gaier80}). To prove orthogonality we will need some properties
of the Zhukowskii function $z=(w+w^{-1})/2$. This function maps the
annulus $1<|w|<\sqrt{\theta}$ one-to-one and conformally onto
$E_{\theta}$ cut along the interval $(-1,1)$ thereby the  circle
$|w|=\sqrt{\theta}$ corresponds to $\partial{E_{\theta}}$. The
inverse function is given by $w=z+\sqrt{z^2-1}$, where the principal
value of the square root is to be chosen. We see by differentiation
that the infinitesimal arc lengths are connected by the relation
$|dz|=|w^2-1||dw|/(2|w|^2)$. Note also that $z^2-1=(w^2-1)^2/(4w^2)$.
Applying the identity $\TS_k(z(w))=w^k+w^{-k}$ we get by the
substitution $z=(w+w^{-1})/2$:
\[
\frac{1}{2\pi}\int\limits_{\partial{E_{\theta}}}\!\!\TS_k(z)\overline{\TS_m(z)}|z^2-1|^{-\frac{1}{2}}|dz|=
\frac{1}{2\pi}\int\limits_{|w|=\sqrt{\theta}}\!\!\left(w^{k}+\frac{1}{w^{k}}\right)
\left(\overline{w}^{m}+\frac{1}{\overline{w}^{m}}\right)\frac{|dw|}{|w|}=
\]\[
=\frac{\theta^{(k+m)/2}}{2\pi}\int\limits_{0}^{2\pi}\!\!e^{i\varphi(k-m)}d\varphi+
\frac{\theta^{-(k+m)/2}}{2\pi}\int\limits_{0}^{2\pi}\!\!e^{i\varphi(m-k)}d\varphi+
\frac{\theta^{(k-m)/2}}{2\pi}\int\limits_{0}^{2\pi}\!\!e^{i\varphi(k+m)}d\varphi+
\]\[
+\frac{\theta^{(m-k)/2}}{2\pi}\int\limits_{0}^{2\pi}\!\!e^{-i\varphi(k+m)}d\varphi
=\left\{\begin{array}{lr} 0, & k\ne m, \\ \theta^{k}+\theta^{-k}, &
k=m\ne 0, \\ 4, & k=m=0.\end{array}\right.
\]
Combined with (\ref{eq:T0defined}) this proves the lemma.$\square$

Denote $\AC={{\cal J}_{\theta}^{\lambda-1/2,\lambda-1/2}}$. We are
ready to formulate our main result for the Gegenbauer expansions.
\begin{theo}\label{theo:AC}
Let $\lambda\geq{0}$.  The space $\AC$ is formed by all restrictions
of the elements of $AL_2(\partial{E_{\theta}})$ to the interval
$(-1,1)$. The norms in $\AC$ and $AL_2(\partial{E_{\theta}})$ are
equivalent.
\end{theo}

\noindent\textbf{Proof.} The proof will be divided in three steps.

\noindent{\em Step} 1. For the space $H_{\GK^{\theta}_0}$ induced by
the kernel (\ref{eq:BilinGigit}) with $\lambda=0$ we want to prove
that
\begin{equation}\label{eq:K0AL2}
H_{\GK^{\theta}_0}=AL_2(\partial{E_\theta}).
\end{equation}
The weight $|z^2-1|^{-1/2}$ is positive and continuous on
$\partial{E_{\theta}}$ so the norms in
$AL_2(\partial{E_{\theta}};|z^2-1|^{-\frac{1}{2}})$ and
$AL_2(\partial{E_{\theta}})$ are equivalent and these spaces coincide
elementwise.  According to Lemma~\ref{lemma:T-ortho} and formula
(\ref{eq:KBilin}) the space
$AL_2(\partial{E_{\theta}};|z^2-1|^{-\frac{1}{2}})$ admits the
reproducing kernel given by
\[
R_\theta(z,\overline{u})=
\frac{1}{4}\sum\limits_{k=0}^{\infty}\frac{\TN_k(z)\overline{\TN_k(u)}}{\theta^{k}+\theta^{-k}}.
\]
This kernel is equivalent to the kernel $\GK^{\theta}_{0}$ due to
(\ref{eq:T0defined}) and inequalities
\[
\frac{1}{2}\sum\limits_{i,j=0}^{n}\GK^{\theta}_{0}(z_i,\overline{z_j})c_i\overline{c_j}\leq
\sum\limits_{i,j=0}^{n}\sum\limits_{k=0}^{\infty}\frac{\TN_k(z_i)\overline{\TN_k(z_j)}}{\theta^{k}+\theta^{-k}}c_i\overline{c_j}\leq
\sum\limits_{i,j=0}^{n}\GK^{\theta}_{0}(z_i,\overline{z_j})c_i\overline{c_j}
\]
satisfied for any choice of $n\in\N$, $c_i\in\C$ and
$z_i\in{E_{\theta}}$.  Hence by Proposition~\ref{prop:K1equalsK2} our
claim is proved.

\noindent{\em Step} 2. Consider the following auxiliary kernel:
\begin{equation}\label{eq:hatKnu}
\hat{K}^{\theta}_{\lambda}(z,\overline{u})
=\frac{\pi^{-1}(\theta^2-1)}{\theta^2-2{\theta}z\overline{u}+1}
{_{2}F_{1}}\left(\frac{\lambda+1}{2},\frac{\lambda+2}{2};\lambda+\frac{1}{2};\frac{4\theta^2(1-z^2)(1-\overline{u}^2)}
{(\theta^2-2{\theta}z\overline{u}+1)^2}\right).
\end{equation}
It is positive definite as will be shown below. Substitution
$\lambda=0$ yields the identity
\begin{equation}\label{eq:hatK0-K0}
\hat{K}^{\theta}_{0}(z,\overline{u})=\GK^{\theta}_{0}(z,\overline{u}).
\end{equation}
We want to prove that for all $\lambda,\mu>-\frac{1}{2}$
\begin{equation}\label{eq:Knu-Kmu}
H_{\hat{K}^{\theta}_{\lambda}}=H_{\hat{K}^{\theta}_{\mu}}.
\end{equation}
According to Proposition~\ref{prop:K1equalsK2} we need to show that
$\hat{K}^{\theta}_{\lambda}\approx\hat{K}^{\theta}_{\mu}$. Following
the definition of the positive definite kernel
(\ref{eq:KPositivity}), choose $n\in\N$, a finite complex sequence
$c_i$ and points $z_i\in E_{\theta}$, $i=\overline{1,n}$. Positive
definiteness of the  kernel $[4\theta^2(1-z^2)(1-\overline{u}^2)]^k/
(\theta^2-2{\theta}z\overline{u}+1)^{2k+1}$ due to its reproducing
property in the Hilbert space of functions representable in the form
$f(z)=(1-z^2)^kg(z)$, where $g$ belongs to the Bergman-Selberg space
generated by the kernel (\ref{eq:B-kernel}) with $\lambda=2k+1$, and
interchange of the order of summations justified by absolute
convergence, lead to the estimates
\begin{eqnarray*}
0\leq\sum\limits_{i,j=1}^{n}\hat{K}^{\theta}_{\mu}(z_i,\overline{z_j})c_i\overline{c_j}=
\frac{\theta^2-1}{\pi}\sum\limits_{k=0}^{\infty}a^{\mu}_k
\sum\limits_{i,j=1}^{n}
\frac{[4\theta^2(1-z_i^2)(1-\overline{z_j}^2)]^k}
{(\theta^2-2{\theta}z_i\overline{z_j}+1)^{2k+1}}c_i\overline{c_j}\hspace{1.5in}
\\
\leq\frac{\theta^2-1}{\pi}\sup_{k\in\N_0}\left\{\frac{a^{\mu}_k}{a^{\lambda}_k}\right\}
\sum\limits_{k=0}^{\infty}a^{\lambda}_k \sum\limits_{i,j=1}^{n}
\frac{[4\theta^2(1-z_i^2)(1-\overline{z_j}^2)]^k}
{(\theta^2-2{\theta}z_i\overline{z_j}+1)^{2k+1}}c_i\overline{c_j}=
\sup_{k\in\N_0}\left\{\frac{a^{\mu}_k}{a^{\lambda}_k}\right\}
\sum\limits_{i,j=1}^{n}\hat{K}^{\theta}_{\lambda}(z_i,\overline{z_j})c_i\overline{c_j},
\end{eqnarray*}
where
\[
a^{\lambda}_k=\frac{([\lambda+1]/2)_k([\lambda+2]/2)_k}{(\lambda+1/2)_k
k!}=
\frac{\Gamma(\lambda+1/2)\Gamma((\lambda+1)/2+k)\Gamma((\lambda+2)/2+k)}
{\Gamma((\lambda+1)/2)\Gamma((\lambda+2)/2)\Gamma(\lambda+1/2+k)k!}>0.
\]
This shows the positive definiteness of the kernel
$\hat{K}^{\theta}_{\mu}$. Using the asymptotic relation
\cite{Bateman53}
\begin{equation}\label{eq:GammaFracAsymp}
\frac{\Gamma(a+z)}{\Gamma(b+z)}=z^{a-b}(1+O(z^{-1})),~~|z|\rightarrow\infty,~~|\arg{z}|<\pi,
\end{equation}
we obtain
\[
\lim_{k\rightarrow\infty}\frac{a^{\mu}_k}{a^{\lambda}_k}=1~\Rightarrow~
0<\sup_{k\in\N_0}\left\{\frac{a^{\mu}_k}{a^{\lambda}_k}\right\}<\infty.
\]
The estimate from below is obtained in the same fashion with
$\sup\limits_{k\in\N_0}\{a^{\mu}_k/a^{\lambda}_k\}$ substituted by
$\inf\limits_{k\in\N_0}\{a^{\mu}_k/a^{\lambda}_k\}$. This proves
equality (\ref{eq:Knu-Kmu}).  Combined with (\ref{eq:hatK0-K0}) and
(\ref{eq:K0AL2}) this gives:
\begin{equation}\label{eq:hatKnuAL2}
H_{\hat{K}^{\theta}_\lambda}=AL_2(\partial{E_\theta})
\end{equation}
for all $\lambda>-1/2$.

\noindent{\em Step} 3. According to (\ref{eq:BilinGigit}) and
(\ref{eq:hatKnu}), the kernel $\GK^{\theta}_{\lambda}$ is related
to the kernel $\hat{K}^{\theta}_{\lambda}$ by
\[
\GK^{\theta}_{\lambda}(z,\overline{u})=B^{\theta}_{\lambda}(z,\overline{u})\hat{K}^{\theta}_{\lambda}(z,\overline{u}),
\]
where
\begin{equation}\label{eq:B-kernel}
B^{\theta}_{\lambda}(z,\overline{u})=\frac{\pi\theta^{2\lambda}}{\tau(\lambda)(\theta^2-2{\theta}z\overline{u}+1)^{\lambda}}=
\frac{h_\lambda^\theta}{(1-\frac{2\theta}{\theta^2+1}z\overline{u})^{\lambda}},
~~h_\lambda^\theta=\frac{\pi\theta^{2\lambda}}{\tau(\lambda)(\theta^2+1)^{\lambda}}.
\end{equation}
For $\lambda>0$ the function $B^{\theta}_{\lambda}(z,\overline{u})$
is the reproducing kernel of the Bergman-Selberg space
$H_{B^{\theta}_{\lambda}}$ \cite{Saitoh97}.  This space comprises
functions holomorphic in the disk $|z|<[(\theta^2+1)/2\theta]^{1/2}$
and having finite norms
\[
\|f\|_{H_{B^{\theta}_{\lambda}}}=\left[\frac{1}{h_\lambda^\theta}
\sum\limits_{k=0}^{\infty}|\mathrm{f}_k|^2\frac{k!}{(\lambda)_k}\left[\frac{\theta^2+1}{2\theta}\right]^k\right]^{\frac{1}{2}},
\]
where $\mathrm{f}_k$ is the $k$-th Taylor coefficient of $f$. The
functions
\[
\gamma_k(z)\bydef[h_\lambda^\theta]^{\frac{1}{2}}\left[\frac{(\lambda)_k}{k!}\right]^{\frac{1}{2}}\left[\frac{2\theta}{\theta^2+1}\right]^{\frac{k}{2}}z^k
\]
constitute a complete orthonormal system in
$H_{B^{\theta}_{\lambda}}$. Note further that the closed ellipse
$\overline{E_\theta}$ is contained in the disk
$|z|<[(\theta^2+1)/2\theta]^{1/2}$ due to inequality
$\sqrt{(\theta^2+1)/(2\theta)}>(\theta^{\frac{1}{2}}+\theta^{-\frac{1}{2}})/2$,
the right hand side of which equals the big semi-axis of the ellipse
$E_\theta$. As stated in Proposition~\ref{prop:TensorDiagonal} the
space $H_{\GK^{\theta}_{\lambda}}$ is obtained by restricting the
elements of the tensor product $H_{B^{\theta}_{\lambda}}\otimes
H_{\hat{K}^{\theta}_\lambda}$ to the diagonal of $E_{\theta}\times
E_{\theta}$ and comprises the functions of the form
\[
f(z)=\sum\limits_{k=1}^{\infty}g_k(z)\gamma_k(z),~~g_k{\in}H_{\hat{K}^{\theta}_\lambda},
~~\sum\limits_{k=1}^{\infty}\|g_k\|_{H_{\hat{K}^{\theta}_\lambda}}^2<\infty.
\]
By (\ref{eq:hatKnuAL2}) we can put $AL_2(\partial{E_\theta})$ instead
of $H_{\hat{K}^{\theta}_\lambda}$ here. For any $g\in
AL_2(\partial{E_\theta})$ consider the estimate
\begin{equation}\label{eq:ggammaEstim}
\|g\gamma_k\|_{AL_2(\partial{E_\theta})}^2=\int\limits_{\partial{E_\theta}}\!\!
|g(z)\gamma_k(z)|^2|dz|\leq
\max_{z\in\partial{E_\theta}}|\gamma_k(z)|^2
\|g\|_{AL_2(\partial{E_\theta})}^2,
\end{equation}
which shows that every product  $g\gamma_k$ belongs to
$AL_2(\partial{E_\theta})$ and hence so does a finite sum of such
products. Denote
\[
\alpha^\theta_{\lambda}(k)\bydef\max_{z\in\partial{E_\theta}}|\gamma_k(z)|=
[h_\lambda^\theta]^{\frac{1}{2}}\left[\frac{(\lambda)_k}{k!}\right]^{\frac{1}{2}}
\left[\frac{2\theta}{\theta^2+1}\right]^{\frac{k}{2}}\left[\frac{\theta+1}{2\sqrt{\theta}}\right]^k=
\]\[
=[h_\lambda^\theta]^{\frac{1}{2}}\left[\frac{(\lambda)_k}{k!}\right]^{\frac{1}{2}}
\left[\frac{\theta+1}{\sqrt{2(\theta^2+1)}}\right]^k,
~~\frac{\theta+1}{\sqrt{2(\theta^2+1)}}<1.
\]
The sequence
\[
S_n(z)=\sum\limits_{k=1}^{n}g_k(z)\gamma_k(z),~~g_k{\in}AL_2(\partial{E_\theta}),
~~\sum\limits_{k=1}^{\infty}\|g_k\|_{AL_2(\partial{E_\theta})}^2<\infty,
\]
is a Cauchy sequence in $AL_2(\partial{E_\theta})$.  Indeed, using
(\ref{eq:ggammaEstim}) and the Cauchy-Schwarz inequality we get
\[
\|S_M-S_N\|_{AL_2(\partial{E_\theta})}^2=
\left\|\sum\limits_{k=N}^{M}g_k\gamma_k\right\|_{AL_2(\partial{E_\theta})}^2=
\sum\limits_{k=N}^{M}\sum\limits_{l=N}^{M}(g_k\gamma_k,g_l\gamma_l)_{AL_2(\partial{E_\theta})}\leq
\]\[
\leq\sum\limits_{k=N}^{M}\sum\limits_{l=N}^{M}\|g_k\gamma_k\|_{AL_2(\partial{E_\theta})}
\|g_l\gamma_l\|_{AL_2(\partial{E_\theta})}\leq
\sum\limits_{k=N}^{M}\sum\limits_{l=N}^{M}\|g_k\|_{AL_2(\partial{E_\theta})}\alpha^\theta_{\lambda}(k)
\|g_l\|_{AL_2(\partial{E_\theta})}\alpha^\theta_{\lambda}(l)=
\]\[
=\left[\sum\limits_{k=N}^{M}\|g_k\|_{AL_2(\partial{E_\theta})}\alpha^\theta_{\lambda}(k)\right]^2\leq
\sum\limits_{k=N}^{M}\|g_k\|_{AL_2(\partial{E_\theta})}^2\sum\limits_{k=N}^{M}[\alpha^\theta_{\lambda}(k)]^2.
\]
Since both $\sum_{k}\|g_k\|_{AL_2(\partial{E_\theta})}^2$ and
$\sum_{k}[\alpha^\theta_{\lambda}(k)]^2$ converge, the above
estimates prove that the sequence $S_n$ is Cauchy.  It follows that
$H_{\GK^\theta_\lambda}\subset AL_2(\partial{E_\theta})$. Inverse
inclusion  $AL_2(\partial{E_\theta}){\subset}H_{\GK^\theta_\lambda}$
is obvious, since $I(z)\equiv{1}$ belongs to
$H_{B^{\theta}_{\lambda}}$ and so for any
$g{\in}AL_2(\partial{E_\theta})$, the product
$Ig=g{\in}H_{\GK^\theta_\lambda}$. $\square$

Now it is not difficult to establish our main result for Jacobi
expansions.
\begin{theo}\label{theo:AJ}
Let $\alpha,\beta\geq-\frac{1}{2}$. Inequality (\ref{eq:sum}) with
$f_k$ defined by (\ref{eq:JacobiCoef}) holds true for all
restrictions to the interval $(-1,1)$ of the elements of
$AL_2(\partial{E_{\theta}})$ and only for them.
\end{theo}
\textbf{Proof.} Choose $\gamma>\max\{\alpha,\beta\}$, then
\[
\JK^\theta_{\alpha,\beta}(z,\overline{u})\ll
M\JK^\theta_{\gamma,\gamma}(z,\overline{u})
\]
for some constant $M>0$.  Indeed, for an arbitrary $n\in\N$, complex
numbers $c_i$ and points $z_i\in E_{\theta}$, $i=\overline{1,n}$,
estimate using (\ref{eq:BilinJacobi}):
\[
\sum\limits_{i,j=1}^{n}\JK^{\theta}_{\alpha,\beta}(z_i,\overline{z_j})c_i\overline{c_j}=
\sum\limits_{k,l=0}^{\infty}a^{\alpha,\beta}_{k,l}\sum\limits_{i,j=1}^{n}
\frac{\theta^{k+l}(1-z_i)^k(1-\overline{z_j})^k(1+z_i)^l(1+\overline{z_j})^l}{(\theta+1)^{2k+2l}}
c_i\overline{c_j}\leq
\]\[
\leq\sup_{k,l\in\N_0}\left\{\frac{a^{\alpha,\beta}_{k,l}}{a^{\gamma,\gamma}_{k,l}}\right\}
\sum\limits_{k,l=0}^{\infty}a^{\gamma,\gamma}_{k,l}\sum\limits_{i,j=1}^{n}
\frac{\theta^{k+l}(1-z_i)^k(1-\overline{z_j})^k(1+z_i)^l(1+\overline{z_j})^l}{(\theta+1)^{2k+2l}}
c_i\overline{c_j}=
\]\[
=\sup_{k,l\in\N_0}\left\{\frac{a^{\alpha,\beta}_{k,l}}{a^{\gamma,\gamma}_{k,l}}\right\}
\sum\limits_{i,j=1}^{n}\JK^{\theta}_{\gamma,\gamma}(z_i,\overline{z_j})c_i\overline{c_j},
\]
where
\[
a^{\alpha,\beta}_{k,l}=\frac{\theta^{\alpha+\beta+1}(\theta-1)}{\tau(\alpha,\beta)
(\theta+1)^{\alpha+\beta+2}}\frac{([\alpha+\beta]/2+1)_{k+l}([\alpha+\beta+3]/2)_{k+l}}
{(\alpha+1)_k(\beta+1)_lk!l!}>0.
\]
Interchange of the order of summations is justified by absolute
convergence of the series (\ref{eq:BilinJacobi}).

Application of formula (\ref{eq:GammaFracAsymp}) yields as
$k,l\rightarrow\infty$
\[
\frac{a^{\alpha,\beta}_{k,l}}{a^{\gamma,\gamma}_{k,l}}=O\left((k+l)^{\alpha+\beta-2\gamma}
k^{\gamma-\alpha}l^{\gamma-\beta}\right)=
O\left((1+l/k)^{\alpha-\gamma}(1+k/l)^{\beta-\gamma}\right)=O(1).
\]
Therefore
$\sup\limits_{k,l\in\N_0}\{a^{\alpha,\beta}_{k,l}/a^{\gamma,\gamma}_{k,l}\}$
is positive and finite. Similarly by choosing
$-1<\eta<\min\{\alpha,\beta\}$ we can prove that
\[
m\JK^{\theta}_{\eta,\eta}\ll \JK^{\theta}_{\alpha,\beta}.
\]
It is left to note that
$\JK^{\theta}_{\beta,\beta}=\GK^{\theta}_{\lambda}$, where
$\lambda=\beta+1/2$, and  $\GK^{\theta}_{\lambda}$ is defined by
(\ref{eq:BilinGigit}).  Now Theorem~\ref{theo:AC} gives the desired
result. $\square$

When $\alpha$ and/or $\beta$ belongs to $(-1,-1/2)$  step 3 of the
proof of the Theorem~\ref{theo:AC} breaks and the problem remains
open.

\begin{corol}\label{corol:AJI}
Condition (\ref{eq:sum}) for the Fourier-Jacobi coefficients
(\ref{eq:JacobiCoef}) of a function $f$ is satisfied for \textbf{all}
$\theta>1$ iff $f$ is the restriction of an entire function to the
interval $(-1,1)$.
\end{corol}

The last theorem and Szeg\"{o}'s theory \cite{Szego91} suggest that
the following much more general conjecture might be true.

\noindent\textbf{Conjecture.}~~{\em Inequality (\ref{eq:sum}) holds
true for the Fourier coefficients in polynomials orthonormal on
$(-1,1)$ with respect to a weight $w$ that satisfies Szeg\"{o}'s
condition $\int_{-1}^{1}\ln{w(x)}dx/\sqrt{1-x^2}>-\infty$ if and only
if $f$ belongs to $AL_2(\partial{E_{\theta}})$. }

\begin{acknowledgments}
The author thanks Professor Martin Muldoon and the York University in
Toronto for hospitality and support during the Fourth ISAAC Congress
in August 2003 and Professor Saburou Saitoh of Gunma University in
Kiryu, Japan, whose book \cite{Saitoh97} on reproducing kernels was
the main inspiration for this research.
\end{acknowledgments}

\end{document}